\documentclass{article}
\usepackage{amsmath,amsfonts,amssymb,theorem}
\mathchardef\ordinarycolon\mathcode`\:     
\mathcode`\:=\string"8000
\def\vcentcolon{\mathrel{\mathop\ordinarycolon}} \begingroup
\catcode`\:=\active \lowercase{\endgroup \let :\vcentcolon }

\newcommand{\F}{{\mathbb F}}

\title{Summary of Delsarte's ``Nombre de Solutions des \'Equations
Polynomiales sur un Corps Fini''}
\author{Wim van Dam\footnote{Massachusetts Institute of Technology,
Center for Theoretical Physics, 77 Massachusetts Avenue,
Cambridge, MA 02139-4307, USA.
email: vandam@mit.edu. 
This work is supported in part by funds provided by the 
U.S.\ Department of Energy  and cooperative research 
agreement DF-FC02-94ER40818 and by a CMI postdoctoral fellowship.
}}
\begin{document}
\maketitle
\begin{abstract}
\noindent 
An English summary is given of Jean Delsarte's article ``Nombre de solutions des \'equations
polynomiales sur un corps fini.'' 
\end{abstract}

\section*{Introduction}
These notes grew out of my desire to check the details of the article: 
\begin{itemize}
\item Jean Delsarte, ``Nombre de solutions des \'equations
polynomiales sur un corps fini'', \emph{S\'eminaire Bourbaki,} 
Expos\'e 39:1--9, March 1951 
\end{itemize}
Because I was only interested in the main result, 
I did not translate the Sections~2 and 4. 
The same notation and equation numbering is maintained and 
the original page numbering is included in the right margin of the text.
Some typos are corrected and I added some potentially helpful comments in italic.
An alternative proof of the result of Delsarte was published 
in \cite{Koblitz}.  This paper cites both Delsarte and 
the 1949 article \cite{Furtado} by Elza Furtado Gomida
as sources for the original result.
I make this translation public to increase the accessibility
of Delsarte's article.
The current translation is by no means authoritative and
does not contain any new results. 
Comments are welcome.

\newpage
\marginpar{page 1}%
\begin{center}
\LARGE{\textbf{Number of Solutions of Polynomial Equations over Finite Fields}}\\
\large{~\\by Jean Delsarte}
\end{center}
\section{Gauss Sums of Finite Fields}
Let $K$ be a finite field $\F_q$ and $K^0$ the multiplicative group $\F_q^\times$.
Also, let $\chi$ be a multiplicative character and $\psi$ a non-trivial
additive character.  The Gauss sum over $K$ is defined by 
\begin{eqnarray*}
g(\chi)&:=&\sum_{x\in K}{\chi(x)\psi(x)}
\end{eqnarray*}
where the $x$ could also run over $K^0$  as $\chi(0)=0$. 
By changing $x$ into $tx$ with $t\in K^0$ we get
\begin{eqnarray*}
g(\chi) &=& \chi(t)\sum_{x\in K}{\chi(x)\psi(tx)},
\end{eqnarray*}
which shows how to convert to different additive characters $\psi$
by multiplying the Gauss sum by a know  factor $\chi(t)$.

A classic results concerns the absolute value of the Gauss sum;
we have 
\begin{eqnarray*}
g(\chi)\bar{g}(\chi) & = & \sum_{x\in K^0}{\sum_{y\in K^0}\chi(xy^{-1})\psi(x-y)}\\
& = & \sum_{x\in K^0}{\chi(x)\sum_{y\in K^0}{\psi((x-1)y)}},
\end{eqnarray*}
where the summation over $y\in K^0$ is $q-1$ if $x=1$, and 
$-1$ otherwise. Thus
\begin{eqnarray*}
|g(\chi)|&=&\sqrt{q}.
\end{eqnarray*}

\section{Finite Extensions of Finite Fields: the Hasse-Davenport Theorem [\dots]}
\marginpar{page 2}%

\marginpar{page 3}%
\section{Some Enumerative Formulae}
Let $E_s$ be the $s$-dimensional $K$ vector space $K^s$, view $E_s$
as a ring with pointwise addition and multiplication. 
Consider the variety defined by the equation
\begin{eqnarray*}
\mathcal{F} & := & \sum_{i=1}^r{a_i x_1^{m_{1i}}\cdots x_s^{m_{si}}} = 0
\end{eqnarray*}
with $a_i \in K^0$ for $i=1,\dots,r$.
The size $q$ of $K$ is big, such that $q-1$ does not divide any of the $m_{ij}$.
Let $\psi$ be an additive character over $K$; we want to 
calculate the summation $S:=\sum_{x}{\psi(\mathcal{F}(x))}$,
where the $x=(x_1,\dots,x_s)\in E_s$. 
We start by calculating the sum $\bar{S}$ where we only sum over those $x$
with $x_j \in K^0$ (that is: $x\in E_s^0$).  Let
\begin{eqnarray}\label{eq1}
y_i & = & x_1^{m_{1i}}\cdots x_s^{m_{si}},
\end{eqnarray}
for $i=1,\dots,r$, where $x=(x_1,\dots,x_s)$ is an invertible 
element of the ring $E_s = K\times\cdots\times K$ (``V\'eron\`ese variety'').
Equation~\ref{eq1} defines a group homomorphism from 
$E_s^0$ to $E_r^0$ (\emph{of the direct products of the multiplicative 
groups $K^0$}).
Let $d$ be the size of the kernel of this homomorphism 
and let $G$ be its image; then we have 
\begin{eqnarray}
\bar{S} = \sum_{x\in E_s^0}{\psi(\mathcal{F}(x))} & = & d\sum_{y\in G}{\psi(a y)}
\end{eqnarray}
with $a=(a_1,\dots,a_r) \in E^0_r$ (\emph{the coefficients of $\mathcal{F}$}).
The product $ay$ is expressed in the ring $E_r$ and the
additive character $\psi$ is extended to this ring
(\emph{according to $\psi(ay):=\psi(a_1y_1+\cdots+a_ry_r)=
\psi(a_1y_1)\cdots\psi(a_ry_r)$}).

Let $\chi$ be a multiplicative character of the group $E_r^0$;
we have for $y=(y_1,\dots,y_r)$ in $E^0_r$ 
\begin{eqnarray*}
\chi(y) = \chi_1(y_1)\dots\chi_r(y_r)
\end{eqnarray*}
\marginpar{page 4}%
where $(\chi_1,\dots,\chi_r)$ is a system of $r$ multiplicative
characters of $K$. 
Let us introduce the group $\tilde{G}$ (orthogonal to $G$),
which is the ensemble of characters on $E_r^0$ with   
$\chi(y)=1$ for every $y\in G$.  Such a character is constant 
on the cosets of $G$ in $E^0_r$.  
(\emph{The group $\tilde{G}$ is the group of multiplicative characters
on $E^0_r/G$; hence its size is $d(q-1)^{r-s}$.})
Now consider the sum
\begin{eqnarray} 
T & = & \sum_{\chi\in\tilde{G}}{\sum_{y\in E^0_r}{\chi(y)\psi(ay)}}.
\end{eqnarray}
For fixed $y$, the sum $\sum_\chi$ is zero when $y$ is outside $G$, 
for $y$ in $G$ the sum is $|\tilde{G}|=d(q-1)^{r-s}$.  Therefore
\begin{eqnarray*}
T & = & d(q-1)^{r-s}\sum_{y\in G}{\psi(ay)},
\end{eqnarray*}
hence
\begin{eqnarray*}
T & = & (q-1)^{r-s}\bar{S},
\end{eqnarray*}
and finally
\begin{eqnarray}
\bar{S} & = & (q-1)^{s-r}\sum_{y\in\tilde{G}}{\sum_{y\in E^0_r}\chi(y)\psi(ay)}. 
\end{eqnarray}
Consider again the Gauss sums, which we have defined for $K$ 
with an additive character $\psi$.  For a multiplicative character 
$\chi=(\chi_1,\dots,\chi_r)$ over $E_r^0$, define   
\begin{eqnarray*}
\mathcal{G}(\chi) & = & g(\chi_1)\cdots g(\chi_r). 
\end{eqnarray*}
Because $\psi(ay) = \psi(a_1y_1)\cdots\psi(a_ry_r)$ we find 
(\emph{because of the earlier derived equality 
$\sum_{y_j}\chi_j(y_j)\psi(a_j y_j) = \bar{\chi}_j(a_j)g(\chi_j)$})
\begin{eqnarray}
\sum_{y\in E_r^0}{\chi(y)\psi(ay)} & = & \bar{\chi}(a)\mathcal{G}(\chi) 
\end{eqnarray}
and hence
\begin{eqnarray}
\bar{S} & = & (q-1)^{s-r}\sum_{\chi\in\tilde{G}}{\bar{\chi}(a)\mathcal{G}(\chi)}.
\end{eqnarray}
\paragraph*{An application of the result.}
Let us try to calculate the number of solutions of $\mathcal{F}=0$ in $E^0_s$, 
which gets denoted by $\bar{N}$. Let
\begin{eqnarray*}
\bar{S}(\psi) & := & \sum_{x\in E_s^0}{\psi(\mathcal{F}(x))}
\end{eqnarray*}
Now calculate the sum of values $\bar{S}(\psi)$ where
$\psi$ ranges over all non-trivial additive characters over $K$:
\marginpar{page 5}%
\begin{eqnarray*}
\sum_\psi{\bar{S}(\psi)} & = & 
\sum_{\psi}{\sum_{x\in E_s^0}{\psi(\mathcal{F}(x))}}. 
\end{eqnarray*}
For fixed $x$, the sum over the characters $\psi$ 
will be $-1$ if $\mathcal{F}(x)$ is not $0$, 
and $q-1$ if $\mathcal{F}(x)=0$, hence
\begin{eqnarray*}
\sum_\psi{\bar{S}(\psi)} & = & (q-1)\bar{N}-((q-1)^s-\bar{N}) \\
& = & q\bar{N}-(q-1)^s
\end{eqnarray*}
because the number values $x\in E_s^0$ for which $\mathcal{F}(x)\neq 0$ 
is $(q-1)^s-\bar{N}$.  Now fix a nontrivial character $\psi_0$. For 
every other nontrivial character $\psi$ we have a $u\in K^0$
such that $\psi(x)=\psi_0(ux)$ for all $x\in K$ (thus establishing a 
bijection between the set of nontrivial characters and $K^0$). 
Moreover we have 
\begin{eqnarray*}
g(\chi_1) = \sum_{x_1\in K^0}{\chi_1(x_1)\psi(x_1)} & = & 
\sum_{x_1\in K^0}{\chi_1(x_1)\psi_0(u x_1)},
\end{eqnarray*}
or 
\begin{eqnarray*}
g(\chi_1) & = & \bar{\chi}_1(u)g_0(\chi_1)
\end{eqnarray*}
where $g_0$ is the Gauss sum where we used the additive character $\psi_0$.
Similarly, one finds 
\begin{eqnarray*}
\mathcal{G}(\chi) & = & \lambda(u)\mathcal{G}_0(\chi),
\end{eqnarray*}
for the Gauss sums over $E_r$. 
(Here $\lambda$ is the multiplicative character over $K^0$ defined by 
the product $\lambda=\chi_1\cdots\chi_r$.)
Finally, one gets
\begin{eqnarray*}
\bar{S}(\psi) & = & 
(q-1)^{s-r}\sum_{\chi\in\tilde{G}}{\lambda(u)\bar{\chi}(a)\mathcal{G}_0(\chi)}.
\end{eqnarray*}
If we want to sum this expression over all nontrivial $\psi$, 
it is sufficient to sum over all $u\in K^0$. 
If $\chi$ is a multiplicative character over $E_r^0$, then 
for the character $\lambda$ over $K^0$, we have that 
$\sum_{u\in K^0}{\lambda(u)}$ is $0$ if $\lambda$ is nontrivial, 
and the sum is $q-1$ if $\lambda$ is trivial. 
We thus have
\begin{eqnarray}
\sum_{\psi}{\bar{S}(\psi)} & = & 
(q-1)^{s-r+1}\sum_{\chi\in\tilde{G}^*}{\chi(a)\mathcal{G}_0(\chi)}
\end{eqnarray}
where $\tilde{G}^*$ is the subgroup of $\tilde{G}$ under the 
restriction that the product $\chi_1\dots\chi_r$ is the 
trivial multiplicative character of $K^0$.
``Without pain'' we thus get the final result
\begin{eqnarray}
\bar{N} & = & 
\frac{1}{q}\left[
(q-1)^s + (q-1)^{s-r+1}\sum_{\chi\in \tilde{G}^*}{\bar{\chi}(a)\mathcal{G}_0(\chi)}
\right].
\end{eqnarray}
\marginpar{page 6}%
\section{The Artin-Weil Series [\dots]}
\marginpar{page 7}%
\marginpar{page 8}%
\marginpar{page 9}%

\end{document}